\documentclass[11pt]{article}
\usepackage{color}
\usepackage{latexsym}
\usepackage{amssymb}
\usepackage{amsmath,amsfonts,theorem,euscript,array,enumerate,mathrsfs}
\usepackage{graphicx}
\newtheorem{Theorem}{Theorem}[part]
\newtheorem{Definition}{Definition}[part]

\newtheorem{Proposition}{Proposition}[part]

\newtheorem{Remark}{Remark}[part]

\def\esssup_#1{\underset{#1}{\mathrm{ess\,sup\, }}}
\def\essinf_#1{\underset{#1}{\mathrm{ess\,inf\, }}}

\def \N{\mathbb{N}}
\def \R{\mathbb{R}}

\def \E{\mathbb{E}}
\def \F{\mathbb{F}}

\def \H{\mathbb{H}}

\def \P{\mathbb{P}}

\def \D{\mathbb{D}}
\def \S{\mathbb{S}}

\def \X{\mathbb{X}}
\def \Y{\mathbb{Y}}

\def \Dc{{\cal D}}

\def \Fc{{\cal F}}

\def \Lc{{\cal L}}

\def \Mc{{\cal M}}

\def \Uc{{\cal U}}

\def \eps{\varepsilon}

\def \ep{\hbox{ }\hfill$\Box$}

\allowdisplaybreaks

\begin{document}

\title{Calculus via regularizations in Banach spaces\\ and Kolmogorov-type path-dependent equations}

\author{
Andrea COSSO\thanks{Laboratoire de Probabilit\'es et Mod\`eles Al\'eatoires, CNRS, UMR 7599, Universit\'e Paris Diderot, France, \sf cosso at math.univ-paris-diderot.fr}
\qquad
Cristina DI GIROLAMI\thanks{Universit\`a di Chieti-Pescara, Dipartimento di
 Economia aziendale, Viale Pindaro 42, I-65127 Pescara, Italy, \sf c.digirolami at unich.it}
\qquad
Francesco RUSSO\thanks{ENSTA ParisTech, Unit\'e de Math\'ematiques appliqu\'ees, 828, boulevard des Mar\'echaux, F-91120 Palaiseau, France, \sf francesco.russo at ensta-paristech.fr}
}
\date{November 28th 2014}

\maketitle

\begin{abstract} The paper reminds the basic ideas of stochastic calculus via regularizations in Banach spaces and its applications to the study of strict solutions of Kolmogorov path dependent equations associated with ``windows'' of
diffusion processes. One makes the link between the Banach space approach
and the so called functional stochastic calculus.
When no strict solutions are available 
one describes the notion of strong-viscosity solution which 
alternative (in infinite dimension) to the classical notion 
of viscosity solution.
\end{abstract}

\vspace{5mm}

\noindent {\bf Key words:} Stochastic calculus via regularization 
in Banach spaces; path dependent Komogorov equation;
functional It\^o calculus.

\vspace{5mm}

\noindent {\bf 2010 Math Subject Classification:}
  60H05;   60H30; 35K10; 35C99.

\newpage

\section{Introduction}

\setcounter{equation}{0} \setcounter{Assumption}{0}
\setcounter{Theorem}{0} \setcounter{Proposition}{0}
\setcounter{Corollary}{0} \setcounter{Lemma}{0}
\setcounter{Definition}{0} \setcounter{Remark}{0}

The present work is a survey (with some new considerations) of recent results on stochastic integration
 in Banach spaces, with applications to Kolmogorov path-dependent partial differential equations (PDEs).

\vspace{1mm}

The extension of It\^o stochastic integration theory for Hilbert valued processes dates only of a few decades, the results of which can be found in the monographs \cite{Metivier82,dpz} and \cite{Walsh} with different techniques. Extension to nuclear valued spaces is simpler and was done in \cite{KaMiW,Ust82}. One of the most natural but difficult situations arises when the processes are Banach space valued.
Big steps forward have been made for instance in \cite{neerven} when the space is of UMD type; on the other hand 
the separable Banach space $C([-T,0])$  of continuous functions $\eta\colon[-T,0]\rightarrow\R$ (endowed with the supremum norm $\|\eta\|_\infty:=\sup_{x\in[-T,0]}|\eta(x)|$) is not UMD.
 This context appears naturally in the study of path-dependent stochastic differential equations (SDEs), as for instance delay equations. An example of such an equation is given by
\begin{equation}
\label{ISDE}
dX_t \ = \ \sigma(t,X_t(\cdot)) dW_t,
\end{equation}
where $W$ is a Brownian motion and $\sigma\colon[0,T] \times C([-T,0]) \rightarrow \R$ 
is continuous and with linear growth. 
Given a continuous real valued process $X$, $X(\cdot)$, also indicated by $\X$, will denote  the so called \emph{window process} associated with $X$, i.e. $X_t(x) := X_{t+x},\,x\in[-T,0]$.
Since $X$ is a continuous process, the natural state space for $\X$ is 
$C([-T,0])$. However, also due to the difficulty of stochastic integration and calculus in that space, most of the authors consider $\X$ as valued in some ad hoc Hilbert space $H$, for example given by the direct sum of $L^2([-T,0])$ and $\R$, see for instance \cite{chojn}.
 To avoid this artificial formulation, a stochastic calculus with $C([-T,0])$-valued stochastic integrators is needed.
 However, if $X = W$ is a classical Brownian motion (therefore we take $\sigma\equiv1$ in \eqref{ISDE}), then the corresponding Brownian window process $\X = X(\cdot)$ has no natural quadratic variation in the sense of Dinculeanu \cite{dincuvisi} or M\'etivier and
Pellaumail \cite{mp}, see Proposition 4.7 in \cite{digirrusso12}.
That quadratic variation is a natural generalization of the one coming from the finite dimensional case.  If $B$ is a separable
Banach space and $\X$ is a $B$-valued process, the notion of quadratic variation (called {\it tensor quadratic variation})
of a process $\X$ introduced by  \cite{dincuvisi} is a process $[\X,\X]$ taking values in the projective tensor
product $B \hat \otimes_\pi B$, see Definition \ref{D:Tensor}. If $B = \R^d$
 and $\X = (X^1, \ldots, X^d)$, $[\X,\X]$  corresponds to the matrix 
$[X^i,X^j]_{1 \le i, j \le d}$.
As mentioned, even though the window Brownian motion does not have a quadratic variation in that sense, it has 
 a more general quadratic variation, known as $\chi$-quadratic variation, first introduced in \cite{DGR} together with the stochastic calculus via regularizations in Banach spaces, for which we also refer to  \cite{DGRnote,digirrusso12,digirfabbrirusso13,DGRNono,DGR2}. The first part of the paper will be devoted to the presentation of the main ideas and results of stochastic calculus via regularizations in Banach spaces, and also to the study of its relation with functional It\^o calculus recently introduced by  \cite{dupire} and \cite{contfournie10, contfournie}.
 
\vspace{1mm}

As an application of this infinite dimensional calculus, we will present a
robust representation of a random variable.  For illustration, let fix 
 $X$ to be a real continuous process  with finite quadratic variation 
$[X]_t=t$, such that $X_0=0$.
Then that representation can be seen as a 
 \emph{robust Clark-Ocone formula}. 
More precisely,
let $h$ be a random variable given by $h=G(\mathbb X_T)$ for some functional $G\colon C([-T,0])\rightarrow\R$. 
We look for a representation (when it is possible) of $h$ of the following type (we remind that $\int_0^T Z_s d^-X_s$ is the forward integral via regularizations defined first in \cite{russovallois93}, which will be recalled in the next section)
\begin{equation}
\label{Repr_h}
G(\mathbb X_T) \ = \ Y_0 + \int_0^T Z_s d^-X_s,
\end{equation}
which, for all $0\leq t\leq T$, can be written as
\begin{equation}
\label{Repr_Y}
Y_t \ = \ G(\mathbb X_T) - \int_t^T Z_s d^-X_s,
\end{equation}
where the pair $(Y,Z)=(Y_t,Z_t)_{t\in[0,T]}$ is required to be adapted to the canonical filtration of $X$. 
The robust aspect is  characterized by the fact that  $Y$ and $Z$ are characterized 
in analytic terms, i.e., through  functions $u,v\colon[0,T]\times C([-T,0])\rightarrow\R$ such that the representation \eqref{Repr_h} becomes
\[
G(\mathbb X_T) \ = \ u(0,\mathbb X_0) + \int_0^T v(s,\mathbb X_s) d^- X_s.
\]
 $u$ and $v$ only depend on the quadratic variation ({\it volatility}) of the process and it
turns out that they are related to the following infinite dimensional partial differential equation:
\begin{equation}
\label{Kolm_Intro}
\begin{cases}
\Lc\,\Uc(t,\eta) \ = \ 0, \;\;\; &\forall\,(t,\eta)\in[0,T[\times C([-T,0]), \\
\Uc(T,\eta) \ = \ G(\eta), &\forall\,\eta\in C([-T,0]),
\end{cases}
\end{equation}
where (we denote by $D_{-t}:=\{(x,x)\colon x\in[-t,0]\}$)
\[
\mathcal{L}\,\Uc(t,\eta) \ = \ \partial_t \Uc(t,\eta)+
\int_{]-t,0]} D^\perp_{dx} \Uc(t,\eta) d^-\eta(x) + \frac{1}{2}\int_{D_{-t}}D^{2}_{dx\,dy} \Uc(t+x,\eta).
\]
Equation \eqref{Kolm_Intro} will be called  Kolmogorov path-dependent PDEs.
This  is the same for all (even non-semimartingale) processes $X$ with the same quadratic variation $[X]_t=t$. As a consequence, this procedure provides a
 Clark-Ocone type representation formula for $h$ which is \emph{robust with respect to quadratic variation}.
In Chapter IV of \cite{mohammed} there is  a characterization of $\Lc$  as
 {\it  infinitesimal generator} (in some weak sense) of the window
 process $\X$, at least for a subspace of the natural subspace that will be considered here. Indeed, the monograph \cite{mohammed} by S.E.A. Mohammed 
constitutes an excellent early
contribution to the theory of functional dependent stochastic differential
equations.

We shall also address the more general problem of characterizing analytically the pair $(Y,Z)$ solution to the following backward stochastic differential equation (here $F\colon[0,T]\times C([-T,0])\times\R\times\R\rightarrow\R$ is a given function)
\[
Y_t \ = \ G(\mathbb X_T) + \int_t^T F(s,\mathbb X_s,Y_s,Z_s) d[X]_s - \int_t^T Z_s d^-X_s,
\]
which is a natural generalization of relation \eqref{Repr_Y}. Another interesting extension corresponds to the case $[X]=\int_0^\cdot \sigma^2(s,\mathbb X_s) ds$, for some function $\sigma\colon[0,T]\times C([-T,0])\rightarrow\R$.

\vspace{1mm}

The last part of the paper is devoted to study more in detail Kolmogorov path-dependent PDEs of the form \eqref{Kolm_Intro} and also of more general type, which naturally arise in stochastic calculus via regularizations in Banach space.
Even in the infinite dimensional case,  Kolmogorov equations is a very active area of research between stochastic calculus and the theory of partial differential equations. On this subject we refer to \cite{cerrai01} and the references therein, and also to \cite{DPZ3} for processes taking values in separable Hilbert spaces, to \cite{fuhrtessi1} for relations with stochastic control theory, to
\cite{flandoligozzi, sobol} for applications to Navier-Stokes equation, and to \cite{flandolidaprato} for connections with infinite dimensional SDEs with irregular drift. Recently, some   interest was devoted to Kolmogorov equations related to Banach space-valued processes, as for instance in \cite{masiero, cerraidaprato}. In the present paper we are interested in Kolmogorov equations on the Banach space $C([-T,0])$, so that the solution is a functional defined on $[0,T]\times C([-T,0])$. $C([-T,0])$ is a natural state space when studying path-dependent stochastic differential equations, as for instance delay equations (even though, as already recalled, the choice of the space $C([-T,0])$ is not usual in the literature, since it is in general more convenient and simpler to work with an Hilbert state space).

We first consider \emph{strict solutions}, namely smooth solutions, to Kolmogorov path-dependent PDEs, for which we 
discuss uniqueness results which are also valid in the case that
$\sigma$ is path-dependent. We recall
 existence  results proved in \cite{cossorusso14} and in \cite{DGRNono}
in the prolongation of
\cite{DGR}.
 Recently, a new approach for existence theorems of smooth solutions has been described in \cite{flandoli_zanco13}. Since, however, strict solutions require quite strong assumptions, we also introduce a weaker notion of solution, called \emph{strong-viscosity solution}, first introduced in \cite{cossorusso14} (we also refer to \cite{cossorusso14_progress} for some new results in this direction), for which we provide a well-posedness result. A strong-viscosity solution 
is defined, in a few words, as the pointwise limit of classical solutions to perturbed equations.
 This definition is similar in spirit to the \emph{vanishing viscosity} method, which represents one of 
the primitive ideas leading to the conception of the modern definition of viscosity solution. This justifies the presence of the term \emph{viscosity} in the name of strong-viscosity solution together with the fact that, as shown in Theorem 3.7 of \cite{cossorusso14}, in the finite dimensional case we have an equivalence result between the notion of strong-viscosity solution and that of viscosity solution.
\vspace{1mm}

The paper is organized as follows. In Section \ref{Regularization_Real_Case} we recall the notion of forward stochastic integral via regularizations for real processes, together with the notion of covariation, and we state the It\^o formula; we end Section \ref{Regularization_Real_Case} with some results on deterministic calculus via regularizations. Section \ref{Regularization_Banach_Case} is devoted to the introduction of stochastic calculus via regularizations in Banach spaces, with a particular attention to the case of window processes; in Section \ref{Regularization_Banach_Case} we also discuss a robust Clark-Ocone formula. Finally, in Section \ref{S:Kolm} we study linear and semilinear Kolmogorov path-dependent equations, we introduce the notions of strict and strong-viscosity solutions, and we investigate their well-posedness.

\section{Stochastic calculus via regularizations}
\label{Regularization_Real_Case}

\setcounter{equation}{0} \setcounter{Assumption}{0}
\setcounter{Theorem}{0} \setcounter{Proposition}{0}
\setcounter{Corollary}{0} \setcounter{Lemma}{0}
\setcounter{Definition}{0} \setcounter{Remark}{0}

\subsection{Generalities}

Let $T\in]0,\infty[$ and consider a probability space $(\Omega,\Fc,\P)$. We denote by $C([-T,0])$ the usual non-reflexive Banach space of continuous functions $\eta\colon[-T,0]\rightarrow\R$ endowed with the supremum norm $\|\eta\|:=\sup_{x\in[-T,0]}|\eta(x)|$. Given a real-valued continuous stochastic process $X=(X_t)_{t\in[0,T]}$ on $(\Omega,\Fc,\P)$, we extend it to all $t\in\R$ as follows: $X_t = X_0$, $\forall\,t<0$, and $X_t=X_T$, $\forall\,t>T$. We then introduce the so-called \textbf{window process} $\X=X(\cdot)$ associated with $X$, which is a $C([-T,0])$-valued stochastic process given by
\[
\X_t := \{X_{t+x},\,x\in[-T,0]\}, \qquad t\in\R.
\]

 Stochastic calculus via regularizations in the finite dimensional
framework
has been largely investigated in the two last decades. It was introduced in \cite{russovallois91, russovallois93} and then developed in several papers (see \cite{russovallois07} for a survey on the subject).
In that calculus, the central object is the forward integral. 
In the present context we will make us of a slightly more general (improper) 
form.

\begin{Definition}
\label{Def_Forw_Int_Real}
Let $X=(X_{t})_{t\in[0,T]}$ and $Y=(Y_{t})_{t\in[0,T]}$ be two real-valued stochastic processes on $(\Omega,\Fc,\P)$, with $X$ continuous and $\int_{0}^{T}|Y_{t}| dt < \infty$ $\P$-a.s.. Suppose that there exists a real continuous process $A=(A_t)_{t \in [0,T]}$ given by 
\begin{equation}\label{IApprox_Real}
A_t \ := \ \lim_{\eps\rightarrow0^+}\int_{0}^{t} Y_s\frac{X_{s+\eps}-X_s}{\eps}ds,	\qquad \forall\,t\in[0,T[,
\end{equation}
where the convergence holds in probability.\begin{enumerate}
\item[\textup{(1)}] The process $A$ will be said  
{\bf forward integral (process) of $Y$ with respect to $X$ (on $[0,T[$)} and it will be denoted by $\int_{0}^\cdot Y d^{-} X$ or $\int_{0}^\cdot Y_s d^{-} X_s$.
\item[\textup{(2)}] If the limit $A_T= \lim_{t \rightarrow T^-} A_t$ holds $\P$-a.s., then $A_T$
will be said 
{\bf (improper) forward integral of $Y$ with respect to $X$ (on $[0,T]$)} and it will be denoted by $\int_{0}^T Y d^{-} X$ or $\int_{0}^T Y_s d^{-} X_s$.
\item[\textup{(3)}] 
For completeness we also remind the {\bf (proper) forward integral of $Y$ with
 respect to $X$ (on $[0,T]$)} as $A_T$ if, in addition to previous two items, 
we have
\[
A_T \ = \ \lim_{\eps\rightarrow0^+}\int_{0}^{T} Y_s\frac{X_{s+\eps}-X_s}{\eps}ds,
\]
where the convergence holds in probability.
\end{enumerate}
\end{Definition}

\begin{Definition} If $I$ is a real subinterval of $[0,T]$,
we say that a family of processes $(H_t^{(\eps)})_{t\in [0,T]}$ converges to $(H_t)_{t\in[0,T]}$ in the \textbf{ucp sense} on $I$,
 if $\sup_{t \in I}|H_t^{(\eps)}-H_t|$ goes to $0$ in probability, as $\eps\rightarrow0^+$. If the interval $I$ will not be mentioned it will be $I = [0,T]$.
\end{Definition}

\begin{Remark} \label{R21}
If the limit \eqref{IApprox_Real} holds in the ucp sense on $[0,T[$ (resp. on
$[0,T]$), then the forward integral $\int_0^\cdot Yd^-X$ of $Y$ with respect to $X$ exists
on $[0,T[$ (resp. ($[0,T]$).
\end{Remark}

We remind now the key notion of covariation. Let us suppose that $Y$, as $X$, is a continuous process.

\begin{Definition}
\label{D33}
The \textbf{covariation of $X$ and $Y$} $($whenever it exists$)$ is given by a 
continuous process $($denoted by $[X,Y]$$)$ such that
\begin{equation} \label{ED33}			
\left[X,Y\right]_{t} \ 
  \ = \ \lim_{\eps\rightarrow 0^{+}}
\frac{1}{\eps} \int_{0}^{t} (X_{s+\eps}-X_{s})(Y_{s+\eps}-Y_{s})ds,
\end{equation}
whenever the limit exists in probability for any $t\in[0,T]$.\\
If $X = Y$, $X$ is called  \textbf{finite quadratic variation process} 
and we denote $[X]:=[X,X]$.
\end{Definition}

If the convergence in \eqref{ED33} holds in the ucp sense
then $[X,Y]$ exists.
We remark that, when $X = Y$, the convergence in probability of \eqref{ED33}
for any $t \in [0,T]$ to $[X,X]$ implies
 that the convergence 
in \eqref{ED33} is also ucp,
 see Lemma 2.1 of \cite{russovallois00}.
 
\vspace{3mm}

Forward integral and covariation are generalizations of the classical It\^o integral and the covariation for semimartingales, as the following result shows (for a proof we refer, e.g., to \cite{russovallois07}). We fix a filtration $\F=(\Fc_t)_{t\in[0,T]}$, $\Fc_T\subset\Fc$, satisfying the usual conditions.

\begin{Proposition}
\label{Pproperties}\quad
\begin{enumerate}
\item[\textup{(i)}] Consider two continuous $\F$-semimartingales $S^1$ and $S^2$. Then, $[S^1,S^2]$ coincides with the standard bracket $[S^1,S^2]=\langle M^1,M^2\rangle$ $($$M^1$ and $M^2$ denote the local martingale parts of $S^1$ and $S^2$, respectively$)$.
\item[\textup{(ii)}] Consider a continuous $\F$-semimartingale $S$ and a c\`{a}dl\`{a}g $\F$-predictable stochastic process $Y$, then the forward integral $\int_0^\cdot Yd^- S$ exists and equals the It\^o integral  $\int_0^\cdot YdS$.
\end{enumerate}
\end{Proposition}

We finally provide It\^o formula in the present finite dimensional setting of stochastic calculus via regularizations, which extends the well-known result for semimartingales to the case of finite quadratic variation processes (see Theorem 2.1 in \cite{russovallois95} for a proof).

\begin{Theorem}[It\^o formula]
\label{ITOFQV}
Let $F\in C^{1,2}\left( [0,T]\times \R;\R \right)$ and consider a real-valued continuous stochastic process $X=(X_t)_{t\in[0,T]}$ with finite quadratic variation. Then,  $\P$-a.s., we have
\begin{align}
\label{E:ITOFQV}
F(t,X_{t}) \ &= \ F(0,X_{0}) + \int_{0}^{t}\partial_t  F(s,X_s)ds + \int_{0}^{t} \partial_{x} F(s,X_s)d^-X_s \\\
&\quad \ + \frac{1}{2}\int_{0}^{t} \partial^2_{x\, x} F(s,X_s)d[X]_s, \notag
\end{align}
for any $0\leq t\leq T$.
\end{Theorem}

\subsection{The deterministic calculus via regularizations}

In the sequel, it will be useful to consider a particular case of finite dimensional stochastic calculus via regularizations, namely the deterministic case which arises when $\Omega$ is a singleton. Let us first fix some useful notations. In this setting we make use of the definite integral on an interval $[a,b]$, where $a<b$ are two real numbers (generally, $a=-T$ or $a=-t$ and $b=0$). We introduce the set $\Mc([a,b])$ of finite signed Borel measures on $[a,b]$. We also denote by $BV([a,b])$ the set of c\`adl\`ag bounded variation functions on $[a,b]$, which is a Banach space when equipped with the norm
\[
\|\eta\|_{BV([a,b])} \ := \ |\eta(b)| + \|\eta\|_{\textup{Var}([a,b])}, \qquad \eta\in BV([a,b]),
\]
where $\|\eta\|_{\textup{Var}([a,b])}=|d\eta|([a,b])$ and $|d\eta|$ is the total variation measure associated to the measure $d\eta\in\Mc([a,b])$ generated by $\eta$: $d\eta(]a,x])=\eta(x)-\eta(a)$, $x\in[a,b]$. Every bounded variation function $f\colon[a,b]\rightarrow\R$ is always suppose to be c\`{a}dl\`{a}g. Moreover, for every function $f\colon[a,b]\rightarrow\R$ we will consider the following two extensions to the entire real line:
\[
f_J(x) :=
\begin{cases}
f(b), & x>b, \\
f(x), & x\in[a,b], \\
f(a), & x<a,
\end{cases} \qquad\qquad
f_{\overline J}(x) :=
\begin{cases}
f(b), & x>b, \\
f(x), & x\in[a,b], \\
0, & x<a,
\end{cases}
\]
where $J:=\,]a,b]$.
\begin{Definition}
\label{D:DeterministicIntegral_closed}
Let $f\colon[a,b]\rightarrow\R$ be a c\`{a}dl\`{a}g function and $g\colon[a,b]\rightarrow\R$ be in $L^1([a,b])$.\\
\textup{(i)} Suppose that the following limit
\[
\int_{[a,b]}g(s)d^-f(s) \ := \ \lim_{\eps\rightarrow0^+}\int_\R g_J(s)\frac{f_{\overline J}(s+\eps)-f_{\overline J}(s)}{\eps}ds,
\]
exists and it is finite. Then, the obtained quantity is denoted by $\int_{[a,b]} gd^-f$ and called \textbf{$($deterministic$)$ forward integral 
of $g$ with respect to $f$ $($on $[a,b]$$)$}.\\
\textup{(ii)} Suppose that the following limit
\[
\int_{[a,b]}g(s)d^+f(s) \ := \ \lim_{\eps\rightarrow0^+}\int_\R g_J(s)\frac{f_{\overline J}(s)-f_{\overline J}(s-\eps)}{\eps}ds,
\]
exists and it is finite. Then, the obtained quantity is denoted by $\int_{[a,b]} gd^+f$ and called \textbf{$($deterministic$)$ backward integral
 of $g$ with respect to $f$ $($on $[a,b]$$)$}.
\end{Definition}

\begin{Definition}
\label{D:DeterministicIntegral_open}
Let $f\colon[a,b]\rightarrow\R$ be a c\`{a}dl\`{a}g function and $g\colon[a,b]\rightarrow\R$ be in $L^1([a,b])$.\\
\textup{(i)} Suppose that the following limit
\[
\int_{]a,b]}g(s)d^-f(s) \ := \ \lim_{\eps\rightarrow0^+}\int_a^b g_J(s)\frac{f_{\overline J}(s+\eps)-f_{\overline J}(s)}{\eps}ds,
\]
exists and it is finite. Then, the obtained quantity is denoted by $\int_{]a,b]} gd^-f$ and called \textbf{$($deterministic$)$ forward integral 
of $g$ with respect to $f$ $($on $]a,b]$$)$}.\\
\textup{(ii)} Suppose that the following limit
\[
\int_{]a,b]}g(s)d^+f(s) \ := \ \lim_{\eps\rightarrow0^+}\int_a^b g_J(s)\frac{f_{\overline J}(s)-f_{\overline J}(s-\eps)}{\eps}ds,
\]
exists and it is finite. Then, the obtained quantity is denoted by $\int_{]a,b]} gd^+f$ and called \textbf{$($deterministic$)$ backward integral
 of $g$ with respect to $f$ $($on $]a,b]$$)$}.
\end{Definition}

Notice that when the two deterministic integrals $\int_{[a,b]}gd^+f$ and $\int_{]a,b]}gd^+f$ exist, they coincide.

\begin{Remark}
{\rm
(i) Let $f\in BV([a,b])$ and $g\colon[a,b]\rightarrow\R$ be a c\`{a}dl\`{a}g function. Then, the forward integral $\int_{]a,b]}gd^-f$ exists and is given by
\[
\int_{]a,b]} g(s) d^-f(s) \ = \ \int_{]a,b]} g(s^-) d\,\!f(s),
\]
where the integral on the right-hand side denotes the classical Lebesgue-Stieltjes integral.\\
(ii) Let $f\in BV([a,b])$ and $g\colon[a,b]\rightarrow\R$ be a c\`{a}dl\`{a}g function. Then, the backward integral $\int_{]a,b]}gd^+f$ exists and is given by
\[
\int_{]a,b]} g(s) d^+f(s) \ = \ \int_{[a,b]} g(s) d\,\!f(s) \ = \ \int_{]a,b]} g(s) d\,\!f(s) + g(a) f(a),
\]
where the integral on the right-hand side denotes the classical 
Lebesgue-Stieltjes integral.
\ep
}
\end{Remark}

\noindent Let us now introduce the deterministic covariation.

\begin{Definition}
\label{D:DeterministicQuadrVar}
Let $f,g\colon[a,b]\rightarrow\R$ be continuous functions and suppose that $0\in[a,b]$. The \textbf{$($deterministic$)$ covariation of $f$ and $g$ $($on $[a,b]$$)$} is defined by
\[						
\left[f,g\right](x) \ = \ \left[g,f\right](x)  \ = \ \lim_{\eps\rightarrow 0^{+}}
\frac{1}{\eps} \int_{0}^x (f(s+\eps)-f(s))(g(s+\eps)-g(s))ds, \qquad x\in[a,b],
\]
if the limit exists and it is finite for every  $x\in[a,b]$. If $f=g$, we set $[f]:=[f,f]$ and it is called \textbf{quadratic variation of $f$ $($on $[a,b]$$)$}.
\end{Definition}
We denote by $V^2$ the set of continuous functions $f\colon[-T,0] \rightarrow \R$
having a deterministic quadratic variation.

Finally, we shall need the following generalization of the deterministic integral when the integrand $g=g(ds)$ is a measure on $[a,b]$ (when the measure $g(ds)$ admits a density with respect to the Lebesgue measure $ds$ on $[a,b]$, we retrieve the deterministic integral introduced in Definition \ref{D:DeterministicIntegral_open}).

\begin{Definition}
\label{D:DeterministicIntegralMeasure}
Let $f\colon[a,b]\rightarrow\R$ be a c\`{a}dl\`{a}g function and $g\in\Mc([a,b])$.\\
\textup{(i)} Suppose that the following limit
\[
\int_{]a,b]}g(ds)d^-f(s) \ := \ \lim_{\eps\rightarrow0^+}\int_{[a,b]} g(ds)\frac{f_{\overline J}(s+\eps)-f_{\overline J}(s)}{\eps},
\]
exists and it is finite. Then, the obtained quantity is denoted by $\int_{]a,b]} gd^-f$ and called \textbf{$($deterministic$)$ forward integral of $g$ with respect to $f$ $($on $]a,b]$$)$}.\\
\textup{(ii)} Suppose that the following limit
\[
\int_{]a,b]}g(ds)d^+f(s) \ := \ \lim_{\eps\rightarrow0^+}\int_{[a,b]} g(ds)\frac{f_{\overline J}(s)-f_{\overline J}(s-\eps)}{\eps},
\]
exists and it is finite. Then, the obtained quantity is denoted by $\int_{]a,b]} gd^+f$ and called \textbf{$($deterministic$)$ backward integral of $g$ with respect to $f$ $($on $]a,b]$$)$}.
\end{Definition}

Indeed, for the sequel, we need to reinforce previous notion.

\begin{Definition}		\label{defn StrExis}
\begin{enumerate}
\item 
We define the following set associated to $\eta\in C([-T,0])$
\begin{equation}	\label{dfn Keta}
K_\eta \ = \ \big\{ \gamma\in C([-T,0])\colon \gamma(x)=\eta(x-\eps),\,x\in [-T,0],\,\eps\in [0,1] \big\}.
\end{equation}
We observe that $K_\eta$ is a compact subset of $C([-T,0])$.\\
\item Let $\Gamma \subset C([-T,0])$.
Let $G\colon[0,T]\times C([-T,0])\rightarrow \mathcal{M}([-T,0])$, $G$ be weakly measurable and bounded. 
We say that 
\begin{equation}
I^-(t,\eta) \ := \ \int_{]-t,0]}G_{dx}(t,\eta)d^-\eta(x) , \;  \quad t\in [0,T],
\end{equation}
$\Gamma$-\textbf{strongly exists} if the following holds for any $ \eta \in \Gamma$.
\begin{itemize}
\item [\textup{(i)}] $\int_{]-t,0]}G_{dx}(t,\eta)d^-\eta(x)$ exists for every $t \in [0,T]$.
\item [\textup{(ii)}]  $K_\eta$ is a subset of $\Gamma$.
For $\eps>0$, $t\in [0,T]$, we set 
$I^-(t,\eta,\eps):=\int_{[-t,0]}G_{dx}(t,\eta)\frac{\eta(x+\eps)-\eta(x)}{\eps}dx$.
We suppose 
that for any $\eta \in \Gamma$, there is $I_\eta\colon[0,T]\rightarrow \R$, 
Lebesgue integrable with respect to $t\in [0,T]$ and such that 
  \begin{equation}       	\label{dfn Upper Bound}
|I^- (t,\gamma,\eps)| \ \leq \ I_\eta(t), \qquad \mbox{for all } \eps\in [0,1], \, t\in [0,T[, \mbox{ and } \gamma \in K_{\eta}.
  \end{equation} 
\end{itemize}
\end{enumerate}
\end{Definition}
Typical choices of $\Gamma$ are the following.
\begin{enumerate}
\item $\Gamma = C([-T,0])$;
\item $\Gamma = V^2$;
\item $\Gamma$ is the linear span 
 of the support of the law of a process $X$.
\end{enumerate}
Sufficient conditions and examples of strong existence of the integrals above are provided in Section 7 of \cite{digirfabbrirusso13}.

We conclude this section by a refinement of the notion
of real finite quadratic variation process.
If $\Gamma = V^2$, a typical example of process $X$
such that $X(\cdot)$ tales values in $\Gamma$ is for instance the
 a $\gamma$-H\"older continuous
process with $\gamma > \frac{1}{2}$, typically a fractional
Brownian motion with Hurst index $ H > \frac{1}{2}$.
If $X$ is a Brownian motion, then $X(\cdot)$ has also a pathwise 
finite quadratic variation,
 see for instance \cite{gradinaru}.
Consequently, if $X$ is the sum of a Wiener process and a H\"older
continuous process with index $\gamma > \frac{1}{2}$, $X(\cdot)$ takes values 
in $V^2$.
 A real process $X$ is said to be of {\bf pathwise  finite
 quadratic variation} if $d\P(\omega)$-a.s. $\eta = X(\omega)$
  belongs to $V^2$
Informally  we can  say that the trajectories of $X$ have a.s.
a 2-variation.

\section{Stochastic calculus via regularizations in Banach spaces}
\label{Regularization_Banach_Case}

\setcounter{equation}{0} \setcounter{Assumption}{0}
\setcounter{Theorem}{0} \setcounter{Proposition}{0}
\setcounter{Corollary}{0} \setcounter{Lemma}{0}
\setcounter{Definition}{0} \setcounter{Remark}{0}

\subsection{General calculus}

In this section we recall briefly basic notions of stochastic calculus for processes $\X$ with values in a Banach space $B$ and its application to window processes $\X=X(\cdot)$, see \cite{DGR,digirrusso12,DGR2} where those notions were introduced. A key ingredient of the stochastic calculus via regularizations in Banach spaces is the notion of Chi-subspace $\chi$, and related $\chi$-covariation. We recall that a Chi-subspace $\chi$ is a (continuously injected) subspace of $(B\hat{\otimes}_\pi B)^{\ast}$, see Definition \ref{D:Chi-subspace} below.

\vspace{3mm}

We begin extending the notion of forward integral introduced in Section \ref{Regularization_Real_Case} for real-valued stochastic processes to the Banach space case. Let $B$ be a separable Banach space equipped with its norm
 $\vert \cdot \vert.$ Given a $B$-valued continuous stochastic process $\mathbb X=(\mathbb X_t)_{t\in[0,T]}$ we extend it to all $t\in\R$ as follows: $\mathbb X_t = \mathbb X_0$, $\forall\,t<0$, and $\mathbb X_t=\mathbb X_T$, $\forall\,t>T$.

\begin{Definition}
\label{def integ fwd}
Consider a $B$-valued stochastic process $\X=(\X_{t})_{t\in[0,T]}$ and a $B^\ast$-valued stochastic process $\Y=(\Y_{t})_{t\in[0,T]}$ on $(\Omega,\Fc,\P)$, with $\X$ continuous and $\int_{0}^{T}\|\Y_{t}\|_{B^{\ast}} dt <\infty$ $\P$-a.s. Suppose that there exists a real continuous process $A=(A_t)_{t \in [0,T]}$ such that 
\begin{equation}\label{IApprox}
A_t \ := \ \lim_{\eps\rightarrow
0^+}\int_{0}^{t} \!\!\textcolor{white}{\Big|_{\textcolor{black}{B^{\ast}}}}\!\Big\langle \Y_s,\frac{\X_{s+\eps} - \X_s}
{\eps}\Big\rangle_{B} ds,	\qquad \forall\,t\in[0,T[,
\end{equation}
where the convergence holds in probability. Then, the process $A$ will be said {\bf forward integral (process) of $\Y$ with respect to $\X$ (on $[0,T[$)} and it will be denoted by $\int_{0}^\cdot {}_{B^{\ast}}{\langle} \Y_s, d^{-} \X_s\rangle_{B}$, or simply by $\int_{0}^\cdot {\langle} \Y_s, d^{-} \X_s\rangle$ when the spaces $B$ and $B^*$ are clear from the context.
 \end{Definition}
		 
When $B=\R$, given a continuous process $X=(X_t)_{t\in[0,T]}$ and a $\P$-a.s. integrable process $Y=(Y_t)_{t\in[0,T]}$, we denote $\int_0^\cdot {}_\R\langle Y,d^-X\rangle_\R$ simply by $\int_0^\cdot Y d^- X$, so we retrieve the forward integral process of $Y$ with respect to $X$ on $[0,T[$ introduced in Definition \ref{Def_Forw_Int_Real}(1).

\vspace{3mm}

Let us now introduce some useful facts about tensor products of Banach spaces.

\begin{Definition}
\label{D:Tensor}
Let $(E,\|\cdot\|_E)$ and $(F,\|\cdot\|_F)$ be two Banach spaces.\\
\textup{(i)} We shall denote by $E\otimes F$ the \textbf{algebraic tensor product} of $E$ and $F$, defined as the set of elements of the form $v = \sum_{i=1}^n e_i\otimes f_i$, for some positive integer $n$, where $e\in E$ and $f\in F$. The map $\otimes\colon E\times F\rightarrow E\otimes F$ is bilinear.\\
\textup{(ii)} We endow $E\otimes F$ with the \textbf{projective norm} $\pi$:
\[
\pi(v) \ := \ \inf\bigg\{\sum_{i=1}^n \|e_i\|_E\|f_i\|_F \ \colon \ v = \sum_{i=1}^n e_i\otimes f_i\bigg\}, \qquad \forall\,v\in E\otimes F.
\]
\textup{(iii)} We denote by $E\hat\otimes_\pi F$ the Banach space obtained as the completion of $E\otimes F$ for the norm $\pi$. We shall refer to $E\hat\otimes_\pi F$ as the \textbf{tensor product of the Banach spaces $E$ and $F$}.\end{Definition}

The definition below was given in \cite{DGR}.

\begin{Definition}
\label{D:Chi-subspace}
Let $E$ be a Banach space. A Banach subspace $(\chi,\|\cdot\|_\chi)$ continuously injected into $(E\hat\otimes_\pi E)^*$, i.e., $\|\cdot\|_\chi\geq\|\cdot\|_{(E\hat\otimes_\pi E)^*}$, will be called a \textbf{Chi-subspace} $($of $(E\hat\otimes_\pi E)^*$$)$.
\end{Definition}

As already mentioned, the notion of Chi-subspace plays a central role in the present Banach space framework, as well as the notion of $\chi$-quadratic variation associated to a Chi-subspace $\chi$, for which we refer to Section 3.2 in \cite{digirrusso12}, and in particular to Definitions 3.8 and 3.9.
If $\X$ is a process admitting $\chi$-quadratic variation, then there exist two maps 
$
[\X]\colon\chi\rightarrow \mathscr{C}([0,T])$ and $
\widetilde{[\X]}\colon\Omega\times [0,T]\rightarrow \chi^{\ast}
$
such that $[\X]$ is linear and continuous, $\widetilde{[\X]}$ has $\P$-a.s. bounded variation and $\widetilde{[\X]}$ is a version of $[\X]$.

\vspace{3mm}

We now present some results of this calculus to window processes, i.e., when $B=C([-T,0])$ and $\X=X(\cdot)$ where $X_t(x) = X_{t+x}$, $\forall\,x \in [-T,0]$. 
A first result about an important integral appearing in the It\^o formula,
in relation with deterministic forward integral via
regularizations, is the following.

\begin{Proposition}		\label{prop 2.3}
Let $\Gamma \subset C([-T,0])$.
Let $\X=X(\cdot)$ be the window process associated with a continuous process $X=(X_t)_{t\in [0,T]}$
such that $\X \in \Gamma$ a.s. 
Let $G$ be weakly bounded and measurable. 
Suppose 
that the forward deterministic integral
\[
I^-(t,\eta) \ := \ \int_{]-t,0]}G_{dx}(t,\eta) d^-\eta(x), \qquad \forall\,
 t \in[0,T],
\]
$\Gamma$-strongly exists. 
Then 
\begin{equation}		\label{eq2.21}
\int_{0}^t \langle G(s,\X_s), d^- \X_{s} \rangle= \int_{0}^t I^-(s,\X_s) ds.
\end{equation}
\end{Proposition}

We will concentrate now  on the Chi-subspace $\chi^0_{Diag}$, which is the following subspace of $C([-T,0])\hat\otimes_\pi C([-T,0])$.
\begin{align*}
\chi^0_{Diag} \ &:= \ \big\{\mu\in\Mc([-T,0]^2)\colon\mu(dx,dy) \ = \ g_1(x,y)dxdy +
\lambda \delta_0(dx)\otimes\delta_0(dy) \\
&\qquad \ + g_2(x)dx\otimes \delta_0(dy) + \delta_0(dx)\otimes g_3(y)dy + g_4(x)\delta_y(dx)\otimes dy, \\
&\qquad \ \ g_1\in L^2([-T,0]^2),\,g_2,g_3\in L^2([-T,0]),\,g_4\in L^\infty([-T,0]),\,\lambda \in\R\big\}.
\end{align*}
 In general, we  refer to the term $g_4(x)\delta_y(dx)\otimes dy$ as the {\bf diagonal component}.

According to Sections 3 and 4 of \cite{digirrusso12}, see also \cite{DGRNono},
one can calculate $\chi$-quadratic variations of a window process
associated with a finite quadratic variation real process.
In particular, we have the following result.

 \begin{Proposition}		\label{prop4A}
Let $X$ be a real finite quadratic variation process
 and $\X=X(\cdot)$ its associated window process. Then $\X=X(\cdot)$ admits a $\chi^0_{Diag}$-quadratic variation which equals $($we denote by $D_{-t}:=\{(x,x)\colon x\in[-t,0]\}$$)$
\begin{equation}	\label{56YH}
\widetilde{[\X]}_t(\mu) \ = \ \mu(\{(0,0)\})[X]_t+\int_{-t}^{0} g_4(x)[X]_{t+x}dx \ = \ \int_{D_{-t}} d\mu(x,y)[X]_{t+x},
\end{equation}
where $\mu$ is a generic element in $\chi^0_{Diag}$ with diagonal component of type $g_4(x)\delta_{y}(dx)dy$, $g_4$ in $L^{\infty}([-T,0])$. 
In particular, if $[X]_t=\int_{0}^t  Z_s  ds$ for an adapted real valued process $(Z_s)_{s\in [0,T]}$, then 
\begin{equation}		\label{56YH2}
\widetilde{[\X]}_t(\mu) \ = \ \int_{0}^t \left(  \int_{D_{-s}} d\mu(x,y)Z_{s+x}\right)ds.
\end{equation}
\end{Proposition}
This allows to state the following theorem, which is an application
 to window processes $\X=X(\cdot)$ of the infinite dimensional It\^o formula
 stated in Theorem 5.2 in \cite{digirrusso12}. In the sequel, $\sigma\colon[0,T]\times C([-T,0])\rightarrow \R$ is a continuous map.
 
\begin{Theorem}  \label{ITO WINDOW}
Let $X$ be a real finite quadratic variation process
 and $\X=X(\cdot)$ its associated window process.
 Let $B=C([-T,0])$ and $F\colon[0,T]\times B\rightarrow \R$ in $C^{1,2}\left( [0,T[\times C([-T,0]) \right)$ in the Fr\'echet sense, such that 
$(t,\eta)\mapsto D^2F(t,\eta)$ is continuous with values in $\chi:=\chi^0_{Diag}$. 
\begin{enumerate}
\item We have
\begin{align}		\label{ItoWindow}
F(t,\X_t) \ &= \ F(0,\X_0)+\int_{0}^t \langle D^\perp_{dx} F(s,\X_s), d^-\X_s \rangle+\int_{0}^t D^{\delta_0} F(s,\X_s)d^-X_s \\
&\quad \ +\frac{1}{2}\int_{0}^{t}\langle D^2F(s,\X_s),d\widetilde{[\X]}_s\rangle, t \in [0,T[,\notag
\end{align}
 whenever either the first or the second integral in the right-hand side exists.
\item If $[X]_t=\int_{0}^t \sigma^2(s,X_s(\cdot)) ds$ then, if  $t \in [0,T[$,
\begin{equation}
\label{Ito_sigma}
\int_{0}^{t}\langle D^2F(s,\X_s),d\widetilde{[\X]}_s\rangle= \int_{0}^t \left( \int_{D_{-s}} D^2 _{dx\, dy} F(s,\X_s)\sigma^2(s+x,\X_{s+x})\right) ds.
\end{equation}
\end{enumerate}
\end{Theorem}
\begin{Remark}
Notice that when the map $F$ in Theorem \ref{ITO WINDOW} satisfies $F(t,\eta)=F(t,\eta(0))$, for all $(t,\eta)\in[0,T]\times C([-T,0])$, so that it does not depend on the ``past'' but only on the ``present value'' of the path $\eta$, then we retrieve It\^o formula \eqref{E:ITOFQV}.
\end{Remark}

\begin{Remark}
As already mentioned, It\^o formula \eqref{ItoWindow} holds if either the first or the second integral in the right-hand side exists. This happens for instance in the two following cases.
\begin{enumerate}
\item $X$ is a semimartingale.
\item $X(\cdot)$ takes values in some subset
$\Gamma$ of $C([-T,0])$ and
 $\int_{]-t,0]}D_{dx}^\perp F(t,\eta)d^- \eta(x)$ $\Gamma$-strongly exists in the sense of Definition \ref{defn StrExis}.
In that case, Proposition \ref{prop 2.3} implies that  
$\int_{0}^{t} {}_{B^{\ast}} {\langle}  
D^\perp F (s,\X_s),  d^- \X_s\rangle_B=\int_0^t I^-(s,\X_s)ds$ as in \eqref{eq2.21}.
\end{enumerate}
\end{Remark}
\textbf{Proof of Theorem \ref{ITO WINDOW}.} \\
Proposition \ref{prop4A} states that $\X$ admits a $\chi^0_{Diag}$-quadratic variation $[\X]$ with version $\widetilde{[\X]}$. 
Item 1. is a consequence of  Theorem 5.2  in \cite{digirrusso12} for $\X=X(\cdot)$. This implies that the forward integral 
$
\int_{0}^{t}
{}_{B^{\ast}} {\langle} DF(s,\mathbb{X}_{s}),d^{-}\mathbb{X}_{s}\rangle_{B}, \ t \in [0,T[,
$
exists and it decomposes into the sum
\begin{equation} \label{ISum}
\int_0^t D^{\delta_0} F(s,\X_s) d^- X_s +   \int_{0}^{t} {}_{B^{\ast}} {\langle}  
D^\perp F (s,\X_s),  d^- \X_s\rangle_B,
\end{equation} provided that at least one of the two addends exists.

Suppose now that $[X]_t=\int_{0}^t \sigma^2(s,X_s(\cdot))ds$. Then  
$\widetilde{[\X]}_{t}(\mu)=\int_{D_{-t}}[X]_{t+x}d\mu(x,y)$ for any $\mu\in \chi^0_{Diag}$.
If $\mu\in \chi^0_{Diag}$, by \eqref{56YH2} setting $Z_s=\sigma^2(s,X_{s}(\cdot))$, we get  
\begin{equation}
\widetilde{[\X]}_t(\mu)=\int_{D_{-t}}\left( \int_{0}^{t+x}
\sigma^2 (s,X_s(\cdot) ) ds\right)d\mu(x,y)=
\int_{0}^t \left( \int_{D_{-s}}d\mu(x,y)\sigma^{2}(s,X_s(\cdot))\right)ds.
\end{equation}
Finally, by elementary integration arguments in Banach spaces 
it follows
\begin{equation}
\int_{0}^t \langle D^2F(s,\X_s), d\widetilde{[\X]}_s\rangle=
\int_{0}^t \left( \int_{D_{-s}} D^2 _{dx\, dy} F(s,\X_s)\sigma^2(s+x,\X_{s+x})\right) ds,
\end{equation}
and the result is established.
\ep

\vspace{3mm}

Now we introduce an important notation.

\begin{Definition}	\label{dfn operator L}
Let $\Uc\colon[0,T]\times C([-T,0])\rightarrow \R$ be in
 $C^{1,2}([0,T[\times C([-T,0]))$. Provided that, for a given $\eta \in C([-T,0])$, $\eta \in C([-T,0)$,
$\int_{]-t,0]} D^\perp_{dx} \Uc(t,\eta) d^-\eta(x)$ $\Gamma$-strongly exists for any  $t \in [0,T[, \eta \in \Gamma$, we define
\begin{align}   \label{operator L}
\mathcal{L}\,\Uc(t,\eta) \ &= \ \partial_t \Uc(t,\eta)+
\int_{]-t,0]} D^\perp_{dx} \Uc(t,\eta) d^-\eta(x) \\
&\quad \ + \frac{1}{2}\int_{D_{-t}}D^{2}_{dx\,dy} \Uc(t+x,\eta)\sigma^2(t+x,\eta(x+\cdot)). \notag
\end{align}

\end{Definition}

\begin{Proposition}	
\label{prop P1}
Let $\Gamma \subset C([-T,0])$. 
Let  $F\colon[0,T]\times C([-T,0])\rightarrow C([-T,0])$ be of class $C^{1,2}\left([0,T[\times C([-T,0]) \right)$ fulfilling the following assumptions.
\begin{itemize}
\item [\textup{(i)}] $\int_{]-t,0]}D^\perp_{dx}F(t,\eta)d^-\eta(x)$, $t\in [0,T[$, $\Gamma$-strongly exists.
\item [\textup{(ii)}] $D^2F\colon[0,T[\times C([-T,0])\rightarrow \chi^0_{Diag}$ exists and it is continuous.
\end{itemize}
Let $X$ be a finite quadratic variation process such that $X(\cdot)$ a.s. lies in $\Gamma$.
\begin{equation}	\label{ECV}
[X]_t=\int_0^{t}\sigma^2(s,\X_s) ds .
\end{equation}
Then,
 the indefinite forward integral $\int_{0}^t D^{\delta_0}F(s,\X_s)d^-X_s, t \in [0,T[$, exists and 
\begin{equation}
\label{eq BN}
F(t,\X_t) \ = \ F(0,\X_0)+\int_{0}^t D^{\delta_0} F(s,\X_s)d^-X_s+\frac{1}{2} \int_{0}^{t} \mathcal{L}F(s,\X_s) ds,
\end{equation}
where $\mathcal{L}F(t,\eta)$ is introduced in Definition \ref{dfn operator L}, see \eqref{operator L}.
\end{Proposition}
{\bf Proof}.
The proof follows from Theorem \ref{ITO WINDOW}, which applies It\^o formula for window processes to $u(s,X_{s}(\cdot))$ between $0$ and $t<T$.
\ep
\vspace{3mm}

Proposition \ref{prop P1}, i.e., the It\^o formula, can be used, in this paper, in two applications. 
\begin{enumerate}
\item To characterize probabilistically the solution of the Kolmogorov equation when $X$ is	a standard stochastic flow.
In particular this is useful to prove uniqueness of strict solutions.
\item To show the robustness representation of a random 
variable, when $X$ is a general finite quadratic variation process.
\end{enumerate}

\subsection{Link with functional It\^o calculus}

Recently a new branch of  stochastic calculus has appeared, known as \emph{functional It\^o calculus}, introduced by  \cite{dupire} and then rigorously developed by 
\cite{contfournie10,contfournie,contfournie13}. It is a stochastic calculus for functionals depending on the all path of a stochastic process, and not only on its current value as in the classical It\^o calculus. One of the main issue of functional It\^o calculus is the definition of the functional (or pathwise or Dupire) derivatives, i.e., the horizontal and vertical derivatives.  
Roughly speaking, the horizontal derivative looks only at the past values of the path, while the vertical derivative looks only at the present value of the path.

In the present section, we shall illustrate how functional It\^o calculus can be interpreted in terms of stochastic calculus via regularizations for window processes. To this end, it will be useful to work within the setting introduced in \cite{cossorusso14}, where functional It\^o calculus was developed by means of stochastic calculus via regularizations. It is worth noting that this is not the only difference between \cite{cossorusso14} and the work \cite{contfournie10} together with \cite{contfournie,contfournie13}. For more information on this point we refer to \cite{cossorusso14}. Here, we just observe that in \cite{contfournie10} it is essential to consider functionals defined on the space of c\`adl\`ag trajectories, since the definition of functional derivatives necessitates of discontinuous paths. Therefore, if a functional is defined only on the space of continuous trajectories (because, e.g., it depends on the paths of a continuous process as Brownian motion), we have to extend it anyway to the space of c\`adl\`ag trajectories, even though, in general, there is no unique way to extend it. In contrast to this approach, in \cite{cossorusso14} it is introduced an intermediate space between the space of continuous trajectories $C([-T,0])$ and the space of c\`adl\`ag trajectories $\D([-T,0])$, denoted $\mathscr C([-T,0])$, which allows to define functional derivatives. $\mathscr C([-T,0])$ is the space of bounded trajectories on $[-T,0]$, continuous on $[-T,0[$ and possibly with a jump at $0$. $\mathscr C([-T,0])$ is endowed with a topology such that $C([-T,0])$ is dense in $\mathscr C([-T,0])$ with respect to this topology. Therefore, any functional $\Uc\colon[0,T]\times C([-T,0])\rightarrow\R$, continuous with respect to the topology of $\mathscr C([-T,0])$, admits a unique extension to $\mathscr C([-T,0])$, denoted $u\colon[0,T]\times\mathscr C([-T,0])\rightarrow\R$. In addition, the time variable and the path have two distinct roles in \cite{cossorusso14}, as for the time variable and the space variable in the classical It\^o calculus. This, in particular, allows to define the horizontal derivative independently of the time derivative, so that, the horizontal derivative defined in \cite{contfournie10} corresponds to the sum of the horizontal derivative and of the time derivative in \cite{cossorusso14}. We mention that an alternative approach to functional derivatives was introduced in \cite{buckdahn_ma_zhang13}.

In the following, we work within the framework introduced in \cite{cossorusso14}. In particular, given a functional $\Uc\colon C([-T,0])\rightarrow\R$ we denote by $D^H\Uc$ and $D^V\Uc$ its horizontal and vertical derivatives, respectively (see Definition 2.11 in \cite{cossorusso14}). Our aim is now to illustrate how the functional derivatives can be expressed in terms of the Fr\'echet derivatives characterizing stochastic calculus via regularizations for window processes. In particular, while it is clear that the vertical derivative $D^V\Uc$ corresponds to $D^{\delta_0}\Uc$, the form of the horizontal derivative $D^H\Uc$ is more difficult to guess. This latter point is clarified by the following two results, which were derived in \cite{cossorusso14}, see Propositions 2.6 and 2.7.

\begin{Proposition}
\label{P:DH=Dacdeta}
Consider a continuously Fr\'{e}chet differentiable map $\Uc\colon C([-T,0])\rightarrow\R$. We make the following assumptions.
\begin{enumerate}
\item[\textup{(i)}] $\forall\,\eta\in C([-T,0])$ there exists $D_x^{\textup{ac}}\Uc(\eta)\in BV([-T,0])$ such that
\[
D_{dx}^\perp \Uc(\eta) \ = \ D_x^{\textup{ac}}\Uc(\eta)dx.
\]
\item[\textup{(ii)}] There exist continuous extensions $($necessarily unique$)$
\[
u\colon\mathscr C([-T,0])\rightarrow\R, \qquad\qquad D_x^{\textup{ac}}u\colon\mathscr C([-T,0])\rightarrow BV([-T,0])
\]
of $\Uc$ and $D_x^{\textup{ac}}\Uc$, respectively.
\end{enumerate}
Then, $\forall\,\eta\in C([-T,0])$,
\begin{equation}
\label{E:DH=Dacdeta}
D^H \Uc(\eta) \ = \ \int_{[-T,0]} D_x^{\textup{ac}} \Uc(\eta) d^+ \eta(x).
\end{equation}
In particular, the horizontal derivative $D^H\Uc(\eta)$ and the backward integral in \eqref{E:DH=Dacdeta} exist.
\end{Proposition}

\begin{Proposition}
\label{P:DH_SecondOrder}
Consider a continuous path $\eta\in C([-T,0])$ with finite quadratic variation on $[-T,0]$. Consider a twice continuously Fr\'echet differentiable map $\Uc\colon C([-T,0])\rightarrow\R$ satisfying
\[
D^2\Uc\colon C([-T,0]) \ \longrightarrow \ \chi_0\subset(C([-T,0])\hat\otimes_\pi C([-T,0]))^*\text{ continuously with respect to $\chi_0$.}
\]
Moreover, assume the following.
\begin{enumerate}
\item[\textup{(i)}] $D_x^{2,Diag}\Uc(\eta)$, the diagonal component of $D^2_x\Uc(\eta)$, has a set of discontinuity which has null measure with respect to $[\eta]$ $($in particular, if it is countable$)$.
\item[\textup{(ii)}] There exist continuous extensions $($necessarily unique$)$$:$
\[
u\colon\mathscr C([-T,0])\rightarrow\R, \qquad\qquad D_{dx\,dy}^2u\colon\mathscr C([-T,0])\rightarrow\chi_0
\]
of $\Uc$ and $D_{dx\,dy}^2\Uc$, respectively.
\item[\textup{(iii)}] The horizontal derivative $D^H\Uc(\eta)$ exists at $\eta$.
\end{enumerate}
Then
\begin{equation}
\label{E:DH=SecondOrder}
D^H \Uc(\eta) \ = \ \int_{]-T,0]} D_{dx}^\perp \Uc(\eta) d^+ \eta(x) - \frac{1}{2}\int_{[-T,0]} D_x^{2,Diag}\Uc(\eta) d[\eta](x).
\end{equation}
In particular, the backward integral in \eqref{E:DH=SecondOrder} exists.
\end{Proposition}

\section{Kolmogorov path-dependent PDE}
\label{S:Kolm}

\setcounter{equation}{0} \setcounter{Assumption}{0}
\setcounter{Theorem}{0} \setcounter{Proposition}{0}
\setcounter{Corollary}{0} \setcounter{Lemma}{0}
\setcounter{Definition}{0} \setcounter{Remark}{0}

\subsection{The framework}

We fix  $\Gamma \subset C([-T,0])$.
Let us consider the following semilinear Kolmogorov path-dependent equation:
\begin{equation}
\label{SemiKolm}
\begin{cases}
\Lc\,\Uc(t,\eta) + F(t,\eta,\Uc, \sigma(t,\eta) D^{\delta_0}\Uc) \ = \ 0, \;\;\; &\forall\,(t,\eta)\in[0,T[\times C([-T,0]), \\
\Uc(T,\eta) \ = \ G(\eta), &\forall\,\eta\in C([-T,0]),
\end{cases}
\end{equation}
where $G\colon C([-T,0])\rightarrow\R$ and $F\colon[0,T]\times C([-T,0])\times\R\times\R\rightarrow\R$ are Borel measurable functions, while the symbol $\mathcal{L}\,\Uc(t,\eta)$ is introduced in Definition \ref{dfn operator L}, see \eqref{operator L}. In the sequel, we think of $\Lc$ as an operator on $C([0,T]\times C([-T,0]))$ with domain
\begin{align*}
\Dc(\Lc) \ := \ \bigg\{\Uc\in C^{1,2}([0,T[\times C([-T,0]))\cap C([0,T]\times C([-T,0]))\colon & \\
\int_{]-t,0]}D_{dx}^\perp \Uc(t,\eta)\,d^-\eta(x) \ \Gamma\text{-strongly  exists} \ \forall\, 
t \in[0,T[ & \bigg\}.
\end{align*}
In the sequel, we will consider the case $\sigma\equiv1$ and give references for more general cases, which are however partly under investigation. 
When $\sigma\equiv1$ we refer to $\Lc$ as \emph{path-dependent heat operator}.

\subsection{Strict solutions}
\label{S42}

We provide the definition of strict solution for equation \eqref{SemiKolm} and we study its well-posedness.

\begin{Definition} \label{dfn strict sol}
We say that $\Uc\colon[0,T]\times C([-T,0])\rightarrow\R$ is a \textbf{strict solution} to the semilinear Kolmogorov path-dependent equation \eqref{SemiKolm} if $\Uc$ belongs to $\Dc(\Lc)$ and solves equation \eqref{SemiKolm}.
\end{Definition}

Concerning the existence and uniqueness of strict solutions, we first consider the linear Kolmogorov path-dependent PDE:
\begin{equation}
\label{Kolm}
\begin{cases}
\Lc\,\Uc(t,\eta) + F(t,\eta) \ = \ 0, \;\;\; &\forall\,(t,\eta)\in[0,T[\times C([-T,0]), \\
\Uc(T,\eta) \ = \ G(\eta), &\forall\,\eta\in C([-T,0]).
\end{cases}
\end{equation}
We have the following uniqueness and existence results for equation \eqref{Kolm}, for which we need to introduce some additional notations. In particular, we consider a complete probability space $(\Omega,\Fc,\P)$ and a real Brownian motion $W=(W_t)_{t\geq0}$ defined on it. We denote by $\F=(\Fc_t)_{t\geq0}$ the natural filtration generated by $W$, completed with the $\P$-null sets of $\Fc$.

\begin{Definition}
\label{E:StochasticFlow}
Let $t\in[0,T]$ and $\eta\in C([-T,0])$. Then, we define the \textbf{stochastic flow}:
\[
\mathbb W_s^{t,\eta}(x) \ = \
\begin{cases}
\eta(x+s-t), &-T \leq x \leq t-s, \\
\eta(0) + W_{x+s} - W_t, \qquad &t-s < x \leq 0,
\end{cases}
\]
for any $t \leq s \leq T$.
\end{Definition}

\begin{Theorem}
\label{T:UniqStrict_Linear}
Let $\Gamma = V^2$.
Consider a strict solution $\Uc$ to \eqref{Kolm} and suppose that there exist two positive constants $C$ and $m$ such that
\begin{equation}
\label{E:PolGrowth_U_Linear}
|G(\eta)| + |F(t,\eta)| + |\Uc(t,\eta)| \ \leq \ C\big(1 + \|\eta\|_\infty^m\big), \qquad \forall\,(t,\eta)\in[0,T]\times C([-T,0]).
\end{equation}
Then, $\Uc$ is given by
\[
\Uc(t,\eta) \ = \ \E\bigg[G(\mathbb W_T^{t,\eta}) + \int_t^T F(s,\mathbb W_s^{t,\eta}) ds\bigg], \qquad \forall\,(t,\eta)\in[0,T]\times C([-T,0]).
\]
In particular, there exists at most one strict solution to the semilinear Kolmogorov path-dependent equation \eqref{SemiKolm} satisfying a polynomial growth
 condition as in \eqref{E:PolGrowth_U_Linear}.
\end{Theorem}
\textbf{Proof.}
Fix $(t,\eta)\in[0,T[\times C([-T,0])$ and $T_0\in[0,T[$. Applying It\^o formula \eqref{ItoWindow} to $\Uc(s,\mathbb W_s^{t,\eta})$ between $t$ and $T_0$, and using \eqref{Ito_sigma}, we obtain
\[
\Uc(t,\eta) \ = \ \Uc(T_0,\mathbb W_{T_0}^{t,\eta}) - \int_t^{T_0} \Lc\,\Uc(s,\mathbb W_s^{t,\eta}) ds - \int_t^{T_0} D^{\delta_0} \Uc(s,\mathbb W_s^{t,\eta}) dW_s. 
\]
Since $\Uc$ solves equation \eqref{Kolm}, we have
\begin{equation}
\label{Uniq_Proof}
\Uc(t,\eta) \ = \ \Uc(T_0,\mathbb W_{T_0}^{t,\eta}) + \int_t^{T_0} F(s,\mathbb W_s^{t,\eta}) ds - \int_t^{T_0} D^{\delta_0} \Uc(s,\mathbb W_s^{t,\eta}) dW_s.
\end{equation}
Consider now the process $M=(M_s)_{s\in[t,T_0]}$ given by
\[
M_s \ := \ \int_t^s D^{\delta_0}\Uc(s,\mathbb W_s^{t,\eta}) dW_s, \qquad \forall\,s\in[t,T_0].
\]
Using the polynomial growth condition of $\Uc$ and $F$, and recalling that, for any $q\geq1$,
\begin{equation}
\label{Uniq_Proof2}
\E\Big[\sup_{t\leq s\leq T}\|\mathbb W_s^{t,\eta}\|_\infty^q\Big] \ < \ \infty,
\end{equation}
we see that $M$ satisfies
\[
\E\Big[\sup_{s\in[t,T_0]}|M_s|\Big] \ < \ \infty.
\]
This implies that $M$ is a martingale. Therefore, taking the expectation in \eqref{Uniq_Proof}, we find
\begin{equation}
\label{UniqLinearProof}
\Uc(t,\eta) \ = \ \E\bigg[\Uc(T_0,\mathbb W_{T_0}^{t,\eta}) + \int_t^{T_0} F(s,\mathbb W_s^{t,\eta}) ds\bigg].
\end{equation}
From the polynomial growth condition \eqref{E:PolGrowth_U_Linear}, together with \eqref{Uniq_Proof2}, we can apply Lebesgue's dominated convergence theorem 
and pass to the limit in \eqref{UniqLinearProof} as $T_0\rightarrow T^-$, from which the claim follows. 
\ep

We remark that previous proof can be easily adapted to the more general
case when $\sigma$ is not necessarily constant.

\begin{Theorem}
\label{T:ExistStrict_Linear}
We suppose $\Gamma = C([-T,0]).$ Let $F\equiv0$
and $G$ admits the cylindrical representation
\begin{equation}
\label{E:G_Cylindrical}
G(\eta) \ = \ g\bigg(\int_{[-T,0]}\varphi_1(x+T)d^-\eta(x),\ldots,\int_{[-T,0]}\varphi_N(x+T)d^-\eta(x)\bigg),
\end{equation}
for some functions $g\in C_p^2(\R^N)$ $($$g$ and its first and second derivatives are continuous and have polynomial growth$)$ and $\varphi_1,\ldots,\varphi_N\in C^2([0,T])$, with $N\in\N\backslash\{0\}$, where the deterministic integrals in \eqref{E:G_Cylindrical} are defined according to Definition \ref{D:DeterministicIntegral_closed}(i). Then, there exists a unique strict solution $\Uc$ to the path-dependent heat equation \eqref{Kolm} satisfying a polynomial growth condition as in \eqref{E:PolGrowth_U_Linear}, which is given by
\[
\Uc(t,\eta) \ = \ \E\big[G(\mathbb W_T^{t,\eta})\big], \qquad \forall\,(t,\eta)\in[0,T]\times C([-T,0]).
\]
\end{Theorem}
\textbf{Proof.}
The proof can be done along the lines of Theorem 3.2 in \cite{cossorusso14}. We simply notice that the idea of the proof is first to show that $\Uc$, as $G$, admits a cylindrical representation. This in turn allows to express $\Uc$ in terms of a function defined on a finite dimensional space: $\Psi\colon[0,T]\times\R^N\rightarrow\R$. Using the regularity of $g$, together with the property of the Gaussian density, we can prove that $\Psi$ is a smooth solution to a certain partial differential equation on $[0,T]\times\R^N$. Finally, using the relation between $\Uc$ and $\Psi$, we conclude that $\Uc$ solves equation \eqref{Kolm}.
\ep

\begin{Remark}
An alternative existence result for strict solutions 
is represented by 
 Proposition 9.53 in \cite{DGR}.
We suppose  \eqref{E:G_Cylindrical} with 
$\varphi_1, \ldots, \varphi_N \in C^2([-T,0])$ such that
\begin{itemize}
\item $g: \R^N \rightarrow \R^N$ in 
only continuous and with linear growth;
\item the matrix $\Sigma_t=(\int_t^T\varphi_i(s)\varphi_j(s)ds)_{1\leq i,j\leq N}$, $\forall\,t\in[0,T]$, has a strictly positive determinant for all $t\in[0,T]$.
\end{itemize}
Then, it follows from Proposition 9.53 in \cite{DGR} that the functional $\Uc$ given by 
\[
\Uc(t,\eta) \ = \ \E\big[G(\mathbb W_T^{t,\eta})\big], \qquad \forall\,(t,\eta)\in[0,T]\times C([-T,0]),
\]
is still the unique strict solution to the path-dependent heat equation \eqref{Kolm} satisfying a polynomial growth condition as in \eqref{E:PolGrowth_U_Linear}.
\end{Remark}
Another existence result is  given below. It is stated and proved in \cite{DGRNono}
and its proof is an adaptation of the proof of
 Theorem 9.41 in \cite{DGR}.


\begin{Theorem}			\label{thm 2 derivate u}
We suppose $\Gamma = C([-T,0]).$
Let $G\in C^{3}\left(C([-T,0])\right)$ such that $D^{3}G$ has polynomial growth. 
Let ${\mathcal U}$ be defined by ${\mathcal U}(t,\eta)=\E \big[ G\big( {\mathbb W}_{T}^{t,\eta} \big)\big]$.
\begin{enumerate}
\item[1)] 
Then $u\in C^{0,2}([0,T]\times C([-T,0]))$. 
\item[2)] Suppose moreover 
\begin{itemize}
\item[i)] $DG(\eta)\in H^{1}([-T,0])$, i.e., function $x\mapsto D_{x}G(\eta)$ is in $H^{1}([-T,0])$, every fixed $\eta$;
\item[ii)]
$DG$ has polynomial growth in $H^{1}([-T,0])$, i.e., 
there is $p\geq 1$ such that 
\begin{equation}			\label{eq HS11}
\eta\mapsto \| DG(\eta) \|_{H^{1}}
\leq 
const\left(  \|\eta\|^{p}_{\infty} +1\right)	\; .		
\end{equation}
\item[iii)] The map 
\begin{equation}		\label{eq HS2}
\eta \mapsto DG(\eta) \hspace{1cm} 
\textrm{considered}  
\hspace{1cm}
C([-T,0])\rightarrow H^{1}([-T,0])
\hspace{1cm} \textrm{is continuous.}
\end{equation} 
\end{itemize}
Then $\mathcal U \in C^{1,2}([0,T]\times C([-T,0]))$ and
 $\mathcal U$ is a strict solution of \eqref{SemiKolm} in the sense of Definition \ref{dfn strict sol}. 
\end{enumerate}
\end{Theorem}
For more existence results concerning strict solutions, with $\sigma$ not necessarily identically equal to 1 and possibly even degenerate,
 we refer to \cite{DGRNono}  and \cite{cossorusso14_progress}.

We end this section proving a uniqueness result for the general semilinear Kolmogorov path-dependent PDE \eqref{SemiKolm}. To this end, we shall rely on the theory of backward stochastic differential equations, for which we need to introduce the following spaces of stochastic processes.

\begin{itemize}
\item $\S^2(t,T)$, $0 \leq t \leq T$, the family of real continuous $\F$-adapted stochastic processes $Y=(Y_s)_{t\leq s\leq T}$ satisfying
\[
\|Y\|_{_{\S^2(t,T)}}^2 := \ \E\Big[ \sup_{t\leq s\leq T} |Y_s|^2 \Big] \ < \ \infty.
\]
\item $\H^2(t,T)$, $0 \leq t \leq T$, the family of  $\R^d$-valued $\F$-predictable stochastic processes $Z=(Z_s)_{t\leq s\leq T}$ satisfying
\[
\|Z\|_{_{\H^2(t,T)}}^2 := \ \E\bigg[\int_t^T |Z_s|^2 ds\bigg] \ < \ \infty.
\]
\end{itemize}

\begin{Theorem}
\label{T:UniqClassical}
Suppose that there exist two positive constants $C$ and $m$ such that
\begin{align*}
|F(t,\eta,y,z) - F(t,\eta,y',z')| \ &\leq \ C\big(|y-y'| + |z-z'|\big), \\
|G(\eta)| + |F(t,\eta,0,0)| \ &\leq \ C\big(1 + \|\eta\|_\infty^m\big),\end{align*}
$\forall\,(t,\eta)\in[0,T]\times C([-T,0])$, $y,y'\in\R$, and $z,z'\in\R$. Consider a strict solution $\Uc$ to \eqref{SemiKolm}, satisfying
\begin{equation}
\label{E:PolGrowth_U}
|\Uc(t,\eta)| \ \leq \ C\big(1 + \|\eta\|_\infty^m\big), \qquad \forall\,(t,\eta)\in[0,T]\times C([-T,0]).
\end{equation}
Then
\[
\Uc(t,\eta) \ = \ Y_t^{t,\eta}, \qquad \forall\,(t,\eta)\in[0,T]\times C([-T,0]),
\]
where $(Y_s^{t,\eta},Z_s^{t,\eta})_{s\in[t,T]} = (\Uc(s,\mathbb W_s^{t,\eta}),D^{\delta_0}\Uc(s,\mathbb W_s^{t,\eta})1_{[t,T[}(s))_{s\in[t,T]}\in\S^2(t,T)\times\H^2(t,T)$ is the solution to the backward stochastic differential equation: $\P$-a.s.,
\[
Y_s^{t,\eta} \ = \ G(\mathbb W_T^{t,\eta}) + \int_s^T F(r,\mathbb W_r^{t,\eta},Y_r^{t,\eta},Z_r^{t,\eta}) dr - \int_s^T Z_r^{t,\eta} dW_r, \qquad t \leq s \leq T.
\]
In particular, there exists at most one strict solution to the semilinear Kolmogorov path-dependent equation \eqref{SemiKolm}.
\end{Theorem}
\textbf{Proof.}
The proof can be done along the lines of Theorem 3.1 in \cite{cossorusso14}, simply observing that the role of the vertical derivative $D^V\Uc$ in \cite{cossorusso14} is now played by $D^{\delta_0}\Uc$.
\ep

\subsection{A robust BSDE representation formula}


 Let $X=(X_t)_{t\in[0,T]}$ be a real process such that its corresponding 
window process 
 $\X = X(\cdot)$ takes
values in $\Gamma = V^2$, i.e. $X$ is a pathwise finite quadratic variation process.
For simplicity we suppose that $[X]_t=t$ and $X_0=0$. 
 Conformally to what we have mentioned in the introduction, given a random variable $h=G(\X_T)$ for some functional $G\colon C([-T,0])\rightarrow\R$, we aim at finding functionals $u,v\colon[0,T]\times C([-T,0])\rightarrow\R$ such that
\[
Y_t \ = \ u(t,\X_t), \qquad Z_t \ = \ v(t,\X_t)
\]
and
\[
Y_t \ = \ G(\mathbb X_T) + \int_t^T F(s,\mathbb X_s,Y_s,Z_s) ds - \int_t^T Z_s d^-X_s,
\]
for all $t\in[0,T]$. In particular, $h$ admits the  representation formula
\[
h \ = \ u(0,\X_0) - \int_0^T F(s,\mathbb X_s,u(s,\X_s),v(s,\X_s)) ds  + \int_0^T v(s,\X_s) d^- X_s.
\]
As a consequence of It\^o formula in Proposition \ref{prop P1}, we have the following result.

\begin{Proposition}
Suppose that $G$ and $F$ are continuous and $u\in C^{1,2}([0,T[\times C([-T,0]))\cap C([0,T]\times C([-T,0]))$. In addition, assume that items (i) and (ii) of Proposition \ref{prop P1} hold with $u$ in place of $F$. Suppose that $u$ solves the Kolmogorov path-dependent PDE \eqref{SemiKolm}. Then 
\begin{equation}
\label{Robust}
h \ = \ Y_0 - \int_0^T F(s,\mathbb X_s,u(s,\X_s),v(s,\X_s)) ds + \int_0^T Z_s d^- X_s,
\end{equation}
with
\[
Y_0 \ = \ u(0,\X_0), \qquad Z_s \ = \ D^{\delta_0}u(s,\X_s).
\]
\end{Proposition}

We refer to \eqref{Robust} as \emph{robust BSDE representation formula} for $h$, and, when $F\equiv0$, as \emph{robust Clark-Ocone formula}.

\subsection{Strong-viscosity solutions}
\label{SubS:StringViscositySolutions}

As we have seen in Section \ref{S42}, we are able to prove an existence result for strict solutions only when the coefficients are regular enough. To deal with more general cases, we need to introduce a weaker notion of solution. We are in particular interested in viscosity-type solutions, i.e., solutions which are not required to be differentiable.

The issue of providing a suitable definition of viscosity solutions for path-dependent PDEs has attracted a great interest. 
We recall that
\cite{ektz},
\cite{etzI,etzII}, and 
 \cite{rtz14} recently provided a definition of viscosity solution to path-dependent PDEs, replacing the classical
 minimum/maximum property, which appears in the standard definition of viscosity solution, with an optimal stopping problem under
 nonlinear expectation~\cite{etzOptStop}. We also recall that other definitions of viscosity solutions for path-dependent PDEs were given by 
\cite{peng12} and 
\cite{tangzhang13}. In contrast with the above cited papers, in the present section we shall adopt the definition of \emph{strong-viscosity solution} introduced in \cite{cossorusso14}, which is not inspired by the standard definition of viscosity solution given in terms of test functions or jets.
Instead, it can be thought, roughly speaking, as the pointwise limit of strict solutions to perturbed equations. We notice that this definition is more similar in spirit to the concept of 
good solution, which turned out to be equivalent to the definition of $L^p$-viscosity solution for certain fully nonlinear partial differential equations, see, e.g., \cite{cerutti_escauriaza_fabes93}, \cite{crandall_kocan_soravia_swiech94}, \cite{jensen96}, and \cite{jensen_kocan_swiech02}. It has also some similarities with the vanishing viscosity method, which represents one of the primitive ideas leading to the conception of the modern definition of viscosity solution. This definition is likewise inspired by the notion of strong solution, as defined for example in \cite{cerrai01}, \cite{gozzi_russo06a}, and \cite{gozzi_russo06b}, even though strong solutions are required to be more regular than strong-viscosity solutions.    We also emphasize that a similar notion of solution, called stochastic weak solution, has been introduced in the recent paper \cite{leao_ohashi_simas14} in the context of variational inequalities for the Snell envelope associated to a non-Markovian continuous process $X$.

A strong-viscosity solution, according to its viscosity nature, is only required to be locally uniformly continuous and with polynomial growth. The term \emph{viscosity} in its name is also justified by the fact that in the finite dimensional case we have an equivalence result between the notion of strong-viscosity solution and that of viscosity solution, see Theorem 3.7 in \cite{cossorusso14}.

We now introduce the notion of strong-viscosity solution for the semilinear Kolmogorov path-dependent equation \eqref{SemiKolm}, which is written in terms of Fr\'echet derivatives, while in \cite{cossorusso14} the concept of strong-viscosity solution was used for an equation written in terms of functional derivatives. Apart from this, the definition we now provide coincides with Definition 3.4 in \cite{cossorusso14}. First, we recall the notion of \emph{locally equicontinuous} collection of functions. 

\begin{Definition}
Let $\mathscr F$ be a collection of $\R^d$-valued functions on $[0,T]\times X$, where $(X,\|\cdot\|)$ is a normed space. We say that $\mathscr F$ is \textbf{locally equicontinuous} if to any $R,\eps>0$ corresponds a $\delta$ such that $|f(t,x)-f(s,y)|<\eps$ for every $f\in\mathscr F$ and for all pair of points $(t,x),(s,y)$ with $|t-s|,\|x-y\|<\delta$ and $\|x\|,\|y\|<R$.
\end{Definition}

\begin{Definition}
\label{D:Strong}
A function $\Uc\colon[0,T]\times C([-T,0])\rightarrow\R$ is called  \textbf{strong-viscosity solution} to the semilinear Kolmogorov path-dependent equation \eqref{SemiKolm} if there exists a sequence $(\Uc_n,G_n,F_n)_n$ satisfying the properties below.
\begin{enumerate}
\item[\textup{(i)}] $\Uc_n\colon[0,T]\times C([-T,0])\rightarrow\R$, $G_n\colon C([-T,0])\rightarrow\R$, and $F_n\colon[0,T]\times C([-T,0])\times\R\times\R\rightarrow\R$ are locally equicontinuous functions such that, for some positive constants $C$ and $m$, independent of $n$,
\begin{align*}
|F_n(t,\eta,y,z) - F_n(t,\eta,y',z')| \ &\leq \ C(|y-y'| + |z-z'|), \\
|\Uc_n(t,\eta)| + |G_n(\eta)| + |F_n(t,\eta,0,0)| \ &\leq \ C\big(1 + \|\eta\|_\infty^m\big),
\end{align*}
for all $(t,\eta)\in[0,T]\times C([-T,0])$, $y,y'\in\R$, and $z,z'\in\R$.
\item[\textup{(ii)}] $\Uc_n$ is a strict solution to
\[
\begin{cases}
\Lc\,\Uc_n \ = \ F_n(t,\eta,\Uc_n,D^{\delta_0}\Uc_n), \;\;\; &\forall\,(t,\eta)\in[0,T)\times C([-T,0]), \\
\Uc_n(T,\eta) \ = \ G_n(\eta), &\forall\,\eta\in C([-T,0]).
\end{cases}
\]
\item[\textup{(iii)}] $(\Uc_n(t,\eta),G_n(\eta),F_n(t,\eta,y,z))\rightarrow(\Uc(t,\eta),G(\eta),F(t,\eta,y,z))$, as $n$ tends to infinity, for any $(t,\eta,y,z)\in[0,T]\times C([-T,0])\times\R\times\R$.
\end{enumerate}
\end{Definition}

\noindent The uniqueness result below for strong-viscosity solution holds.

\begin{Theorem}
Let $\Uc\colon[0,T]\times C([-T,0])\rightarrow\R$ be a strong-viscosity solution to the semilinear Kolmogorov path-dependent equation \eqref{SemiKolm}. Then
\[
\Uc(t,\eta) \ = \ Y_t^{t,\eta}, \qquad \forall\,(t,\eta)\in[0,T]\times C([-T,0]),
\]
where $(Y_s^{t,\eta},Z_s^{t,\eta})_{s\in[t,T]}\in\S^2(t,T)\times\H^2(t,T)$, with $Y_s^{t,\eta}=\Uc(s,\mathbb W_s^{t,\eta})$, solves the backward stochastic differential equation: $\P$-a.s.,
\[
Y_s^{t,\eta} \ = \ G(\mathbb W_T^{t,\eta}) + \int_s^T F(r,\mathbb W_r^{t,\eta},Y_r^{t,\eta},Z_r^{t,\eta}) dr - \int_s^T Z_r^{t,\eta} dW_r, \qquad t \leq s \leq T.
\]
In particular, there exists at most one strong-viscosity solution to the semilinear Kolmogorov path-dependent equation \eqref{SemiKolm}.
\end{Theorem}
\textbf{Proof.}
Let us give only a sketch of the proof (for a similar argument and more details, see Theorem 3.3 in \cite{cossorusso14}). Consider a sequence $(\Uc_n,G_n,F_n)_n$ satisfying conditions (i)-(ii)-(iii) of Definition \ref{D:Strong}. For every $n\in\N$ and any $(t,\eta)\in[0,T]\times C([-T,0])$, we know from Theorem \ref{T:UniqClassical} that $(Y_s^{n,t,\eta},Z_s^{n,t,\eta})_{s\in[t,T]} = (\Uc_n(s,\mathbb W_s^{t,\eta}),D^{\delta_0}\Uc_n(s,\mathbb W_s^{t,\eta}))_{s\in[t,T]}\in\S^2(t,T)\times\H^2(t,T)$ is the solution to the backward stochastic differential equation: $\P$-a.s.,
\[
Y_s^{n,t,\eta} \ = \ G_n(\mathbb W_T^{t,\eta}) + \int_s^T F_n(r,\mathbb W_r^{t,\eta},Y_r^{n,t,\eta},Z_r^{n,t,\eta}) dr - \int_s^T Z_r^{n,t,\eta} dW_r, \qquad t \leq s \leq T.
\]
Thanks to a limit theorem for BSDEs (see Proposition C.1 in \cite{cossorusso14}), and using the hypotheses on the coefficients, we can pass to the limit in the above backward equation as $n\rightarrow\infty$, from which the thesis follows.
\ep

\vspace{3mm}

We finally address the existence problem for strong-viscosity solutions in the linear case, and in particular when $F\equiv0$.

\begin{Theorem}
\label{T:Exist}
Let $F\equiv0$ and $G\colon C([-T,0])\rightarrow\R$ be a locally uniformly continuous map satisfying
\[
|G(\eta)| \ \leq \ C(1 + \|\eta\|_\infty^m), \qquad \forall\,\eta\in C([-T,0]),
\]
for some positive constants $C$ and $m$. Then, there exists a unique strong-viscosity solution $\Uc$ to equation \eqref{SemiKolm}, which is given by
\[
\Uc(t,\eta) \ = \ \E\big[G(\mathbb W_T^{t,\eta})\big], \qquad \forall\,(t,\eta)\in[0,T]\times C([-T,0]).
\]
\end{Theorem}
\textbf{Proof.}
The proof can be done along the lines of Theorem 3.4 in \cite{cossorusso14}. Let us give an idea of it. We first fix $\eta\in C([-T,0])$ and derive a Fourier series expansion of $\eta$ in terms of a smooth orthonormal basis of $L^2([-T,0])$. This allows us to approximate $G$ with a sequence of functions $(G_n)_n$, where $G_n$ depends only on the first $n$ terms of the Fourier expansion of $\eta$. Noting that the Fourier coefficients can be written in terms of a forward integral with respect to $\eta$, we see that every $G_n$ has a cylindrical form. Moreover, even if $G_n$ is not necessarily smooth, we can regularize it. After this final smoothing, we end up with a terminal condition, that we still denote $G_n$, which is smooth and cylindrical. As a consequence, from Theorem \ref{T:ExistStrict_Linear} it follows that the corresponding Kolmogorov path-dependent equation admits a unique strict solution $\Uc_n$ given by
\[
\Uc_n(t,\eta) \ = \ \E\big[G_n(\mathbb W_T^{t,\eta})\big], \qquad \forall\,(t,\eta)\in[0,T]\times C([-T,0]).
\]
It is then easy to show that the sequence $(\Uc_n,G_n)_n$ satisfies points (i)-(ii)-(iii) of Definition \ref{D:Strong}, from which the thesis follows.
\ep

\medskip

{\bf ACKNOWLEDGEMENTS}. The third named author benefited from the
support of the ``FMJH Program Gaspard Monge in optimization and operation
research'' (Project 2014-1607H) and from the support to this program from
EDF. The second name author was partially supported by the Fernard
 Braudel-IFER outgoing fellowship, funded by the Fondation Maison de Science 
de l'Homme and the European Commission, Action Marie Curie COFUND, 7e PCRD.

\small
\bibliographystyle{plain}
\bibliography{biblio}

\end{document}